\documentclass{amsproc}
\usepackage{amsmath,amssymb,latexsym,amsfonts,cite,hyperref}
\usepackage[all]{xy}

\hypersetup{unicode=false, pdftoolbar=true, pdfmenubar=true, pdffitwindow=false, pdfstartview={FitH}, pdfnewwindow=true, colorlinks=false, linkcolor=blue, citecolor=green}

\newcommand{\X}{\mathbb X}

\renewcommand{\O}{\mathcal O}
\newcommand{\A}{\mathcal A}

\newtheorem{theorem}{Theorem}[section]
\newtheorem{proposition}[theorem]{Proposition}

\newtheorem{corollary}[theorem]{Corollary}

\newtheorem{lemma}[theorem]{Lemma}
\numberwithin{equation}{section}
\theoremstyle{definition}

\begin{document}
\title[FRT-Construction via Quantum Affine Algebras and Smash Products]{The FRT-Construction via Quantum Affine Algebras and Smash Products}
\author{Garrett Johnson and Chris Nowlin}
\address{Department of Mathematics, University of California, Santa Barbara, CA, 93106}
\email{johnson@math.ucsb.edu\\ cnowlin@math.ucsb.edu}

\begin{abstract}
\noindent For every element $w$ in the Weyl group of a simple Lie algebra $\mathfrak{g}$, De Concini, Kac, and Procesi defined a subalgebra ${\mathcal U}_q^w$ of the quantized universal enveloping algebra ${\mathcal U}_q(\mathfrak{g})$. The algebra ${\mathcal U}_q^w$ is a deformation of the universal enveloping algebra ${\mathcal U}(\mathfrak{n}_+\cap w.\mathfrak{n}_-)$. We construct smash products of certain finite-type De Concini-Kac-Procesi algebras to obtain ones of affine type; we have analogous constructions in types $A_{n}$ and $D_n$. We show that the multiplication in the affine type De Concini-Kac-Procesi algebras arising from this smash product construction can be twisted by a cocycle to produce certain subalgebras related to the corresponding Faddeev-Reshetikhin-Takhtajan bialgebras.
 \end{abstract}

\maketitle

\section{Introduction} Let $k$ be an infinite field and suppose an algebraic $k$-torus $H$ acts rationally on a noetherian $k$-algebra $A$ by $k$-algebra automorphisms.  Goodearl and Letzter \cite{GoodLet} showed that $\text{spec}(A)$ is partitioned into strata indexed by the $H$-invariant prime ideals of $A$. Furthermore, they showed that each stratum is homeomorphic to the prime spectrum of a Laurent polynomial ring.  The Goodearl-Letzter stratification results apply to the case when $A$ is an iterated Ore extension under some assumptions relating the action of $H$ to the structure of $A$.  In this setting Cauchon's deleting derivations algorithm \cite{C} gives an iterative procedure for classifying the $H$-primes. After several such algebras were studied, such as the algebras of quantum matrices $\O_q(M_{\ell ,p}(k))$ \cite{C, GL, La}, it was noticed that many of these algebras fall into the setting of De Concini-Kac-Procesi algebras \cite{DKP}.

The  De Concini-Kac-Procesi algebras are subalgebras of quantized universal enveloping algebras ${\mathcal U}_q(\mathfrak{g})$ associated to the elements of the corresponding Weyl group $W_\mathfrak{g}$.  They may be viewed as deformations of the universal enveloping algebra ${\mathcal U}(\mathfrak{n}_+\cap w.\mathfrak{n}_-)$, where $\mathfrak{n}_+$ and $\mathfrak{n}_-$ are the positive and negative nilpotent Lie subalgebras of ${\mathfrak g}$, respectively.  M\'eriaux and Cauchon \cite{MC} and Yakimov \cite{Y} recently proved that the poset of $H$-primes of a De Concini-Kac-Procesi algebra ${\mathcal U}_q^w$ ordered under inclusion is isomorphic to the poset $W^{\leq w}$ of Weyl group elements less than or equal to $w$ under the Bruhat ordering. In \cite{Y}, Yakimov also gives explicit generating sets for the $H$-primes in terms of Demazure modules.

In this paper we introduce three type-$D$ algebras. The first algebra is obtained from a De Concini-Kac-Procesi algebra of type $\mathfrak{so}_{2n+2}$ and a smash product construction. We show that this algebra is isomorphic to our second algebra, a De Concini-Kac-Procesi algebra associated to the affine Weyl group of type $\widehat{D}_{n+1}$. Finally, we show that twisting the multiplication in these algebras by a certain $2$-cocycle produces algebras related to the type-$D$ Faddeev-Reshetikhin-Takhtajan bialgebras. In the last section, we produce analogous results with type-$A$ algebras.  In a forthcoming publication we will return to the $H$-spectrum of these algebras.

Section \ref{section smashD} introduces the first of these algebras, an algebra which resembles a smash product of a De Concini-Kac-Procesi algebra with itself.   Let $k$ be an algebraically closed field of characteristic zero and suppose $q\in k$ is not a root of unity.  Fix an integer $n\geq 3$. Let $W(D_{n+1})$ be the Weyl group of type $D_{n+1}$ with standard generating set $\left\{ s_1,\dots,s_{n+1}\right\}$ and let \begin{equation}w_n=(s_{n+1}s_n\cdots s_2s_1)(s_3s_4\cdots s_ns_{n+1})\in W(D_{n+1}). \end{equation} Let ${\mathcal U}_D^{\geq 0}$ denote the quantized positive Borel algebra of type $D_{n+1}$ and let ${\mathcal U}_q^{w_n}$ be the De Concini-Kac-Procesi subalgebra of ${\mathcal U}_D^{\geq 0}$ corresponding to $w_n$.  In fact,  ${\mathcal U}_q^{w_n}$ is isomorphic to $\O_q(\mathfrak{o}k^{2n})$, the algebra of even-dimensional quantum Euclidean space.  We define an action $\lambda$ of ${\mathcal U}_D^{\geq 0}$ on ${\mathcal U}_q^{w_n}$, which is a modification of the adjoint action of the Hopf algebra ${\mathcal U}_D^{\geq 0}$ on itself. This action equips ${\mathcal U}_q^{w_n}$ with the structure of a left ${\mathcal U}_D^{\geq 0}$-module algebra.  We then consider the smash product ${\mathcal U}_q^{w_n}\# {\mathcal U}_D^{\geq 0}$ with respect to $\lambda$ and set $({\mathcal U}_q^{w_n})^\#$ to be the subalgebra of  ${\mathcal U}_q^{w_n}\# {\mathcal U}_D^{\geq 0}$ generated by $\{ u\#1,\ 1\#u\mid u\in {\mathcal U}_q^{w_n}\}$.  This is the first of the three type-$D$ algebras.

In Section \ref{section affD} we introduce a second type-$D$ algebra; it is a De Concini-Kac-Procesi algebra of affine type.  Let $ W(\widehat D_{n+1})$ denote the affine Weyl group of type $\widehat D_{n+1}$ with generating set $\{s_0,s_1,\dots,s_{n+1}\}$ and let \begin{equation}\widehat w_n=(s_{n+1}\cdots s_1)(s_3\cdots s_{n+1})s_0(s_n\cdots s_3)(s_1\cdots s_n)s_0\in W(\widehat D_{n+1}).\end{equation} The main result of this section is Theorem \ref{affineD=smash}, where we prove that the algebras ${\mathcal U}_q^{\widehat w_n}$ and $({\mathcal U}_q^{w_n})^\#$ are isomorphic.

In Section \ref{section X algebra} we introduce an algebra which we denote by $\mathbb{X}_{n,q}$. We show that $\mathbb{X}_{n,q}$ is related to the bialgebra $\mathcal A(R_{D_n})$ arising from the type-$D_n$ FRT construction.  In particular, we label the standard generators of  $\mathcal A(R_{D_n})$ by $Y_{ij}$, for $1\leq i,j\leq 2n$, and let $T_{2,n}\subseteq \mathcal A(R_{D_n})$ be the subalgebra generated by $\{Y_{ij}:1\leq i\leq 2,\ 1\leq j\leq 2n\}$ and observe that there is a surjective algebra homomorphism $\mathbb{X}_{n,q}\to T_{2,n}$ (see Proposition \ref{surjective homomorphism}).  We thus refer to $\mathbb{X}_{n,q}$ as a \emph{parent} of $T_{2,n}$.  Finally, in Thm. \ref{twisting} we prove that $\mathbb{X}_{n,q}$ is isomorphic to a cocycle twist (in the sense of \cite{AST91}) of ${\mathcal U}_q^{\widehat{w}_n}$.  From this, it follows that $\mathbb{X}_{n,q}$ is an iterated Ore extension over $k$.

In Section \ref{section A analogue} we proceed to demonstrate analogous results in the type $A_m$ setting. We fix an integer $m>1$ and  let $W(A_m)$ be the Weyl group of type $A_m$ with generating set $\left\{ s_1,\dots,s_m\right\}$.  Let \begin{equation} c_m=s_1\cdots s_m\in W(A_m) \end{equation} denote a Coxeter element.  Notice that the De Concini-Kac-Procesi algebra ${\mathcal U}_q^{c_m}$ is isomorphic to ${\mathcal O}_q(k^m)$, the algebra known as quantum affine space.  Quantum euclidean space, seen in Section 3, can be thought of as a type-$D$ analogue of quantum affine space.  Let ${\mathcal U}_A^{\geq 0}$ denote the quantized positive Borel algebra of type $A_m$.  We define an action $\lambda_A:{\mathcal U}_A^{\geq 0}\otimes {\mathcal U}_q^{c_m}\to{\mathcal U}_q^{c_m}$ endowing ${\mathcal U}_q^{c_m}$ with the structure of a left ${\mathcal U}_A^{\geq 0}$-module algebra and define $({\mathcal U}_q^{c_m})^\#$ to be the subalgebra of ${\mathcal U}_q^{c_m}\#{\mathcal U}_A^{\geq 0}$ generated by $\{1\# u,\ u\# 1\mid u\in {\mathcal U}_q^{c_m}\}.$ Finally, we let $W(\widehat A_m)$ denote the affine Weyl group of type $\widehat A_m$ with generating set $\{s_0,s_1,\dots,s_{m}\}$  and let \begin{equation}\widehat c_m=(s_1\cdots s_m)(s_0s_1\cdots s_{m-1})\in  W(\widehat A_m).\end{equation} In Thm. \ref{affineA=smash}, we prove that the corresponding De Concini-Kac-Procesi algebra ${\mathcal U}_q^{\widehat{c}_m}$ is isomorphic to $({\mathcal U}_q^{c_m})^\#$. We further demonstrate that ${\mathcal U}_q^{\widehat c_m}$ is isomorphic to a cocycle twist of ${\mathcal O}_q(M_{2,m})$, the algebra of $2\times m$ quantum matrices.  We see $\X_{m,q}$ as a ``type D'' analogue of ${\mathcal O}_q(M_{2,m})$ because ${\mathcal O}_q(M_{2,m})$ is a subalgebra of the FRT-bialgebra $\mathcal A(R_{A_{m-1}})\cong{\mathcal O}_q(M_m(k))$. The key distinction is that $\O_q(M_{2,m}(k))$ is a subalgebra of $\A(R_{A_{m-1}})$, whereas $\X_{m,q}$ is a parent of the analogous subalgebra $T_{2,m}\subseteq\A(R_{D_m})$.

\vspace{.5 cm}

\noindent{\bf Acknowledgements.}
We would like to acknowledge Milen Yakimov and Ken Goodearl for helpful conversations.

\section{Preliminaries}

Let $Q$ be a $\mathbb{Z}$-module with basis $\Pi=\{\alpha_1,...,\alpha_n\}$.   Suppose $\left<\,,\,\right>$ is a symmetric form on $Q$  with \begin{equation} c_{ij}:=  2\left<\alpha_i,\alpha_j\right>/\left<\alpha_i,\alpha_i\right>\end{equation} non-positive integers for every $i,j\in\{1,...,n\}$ with $i\neq j$, and there exist coprime positive integers $d_1,...,d_n$ so that the matrix $(d_ic_{ij})$ is symmetric. Let $k$ be a field and let $q\in k$ be nonzero. Set $q_i=q^{d_i}$ and assume $q_i\neq \pm 1$. To the triple $(Q,\Pi, q)$, we have an associated quantized enveloping algebra ${\mathcal U}$. As a $k$-algebra, ${\mathcal U}$ is generated by $u_1^{\pm},...,u_n^{\pm}$ and $\{v_{\mu}: \mu\in Q\}$ and has defining relations
\begin{align}
 & v_0 = 1, \hspace{1cm} v_\mu v_\rho=v_{\mu+\rho},    &  (\mu, \rho\in Q),& \\
 & v_\mu u_i^{\pm} = q^{\pm \left< \mu ,\alpha_i\right>} u_i^{\pm}v_\mu,        &  (\mu\in Q, i\in\{1,...,n\}),& \\
 &u_i^+u_j^- = u_j^-u_i^++\delta_{ij}\frac{v_{\alpha_i}-v_{-\alpha_i}}{q_i-q_i^{-1}},    &  (i,j\in\{1,...,n\}),& \\
 & \sum_{r=0}^{1-c_{ij}} (-1)^r\left[\begin{array}{c}1-c_{ij} \\ r\end{array}\right]_{q_i}(u_i^\pm)^{1-c_{ij}-r} u_j^\pm(u_i^\pm)^r =0, & (i\neq j).
\end{align}

Here, \begin{equation}[\ell ] _{q_i} = \frac{q_i^\ell -q_i^{-\ell}}{q_i-q_i^{-1}}, \hspace{.6cm}[\ell]_{q_i}! = [1]_{q_i}\cdots [\ell]_{q_i},\hspace{.6cm}\left[\begin{array}{c}\ell \\ m\end{array}\right]_{q_i}=\frac{[\ell]_{q_i}!}{[m]_{q_i}![\ell - m]_{q_i}!}.\end{equation}

Furthermore, ${\mathcal U}$  has a Hopf algebra structure with comultiplication $\Delta$, antipode $S$, and counit $\epsilon$ maps given by
\begin{align}
&\label{comultiplication formula}\Delta (u_i^+) = v_{-\alpha_i}\!\otimes\! u_i^+\!+\!u_i^+\!\otimes \!1,&                    &\hspace{-.1cm}\Delta (v_\mu ) = v_\mu\!\otimes \!v_\mu, &             &\hspace{-.1cm}\Delta(u_i^-)= 1\!\otimes\! u_i^-\!+\!u_i^-\!\otimes \!v_{\alpha_i},\\
&S (u_i^+) = -v_{\alpha_i}E_i,&                                                                  &S(v_\mu )=v_{-\mu},&                                               &S (u_i) = -u_i^-v_{-\alpha_i},\\
&\epsilon (u_i^+) = 0,&                                                                                  &\epsilon (v_\mu ) =1,&                                              &\epsilon (u_i^-)=0,\end{align}
for every $\mu\in Q$ and $1\leq i\leq n$.

For every $i\in\{1,...,n\}$, we let $s_i:Q\to Q$ be the simple reflection \begin{equation}s_i: \mu\to \mu-\frac{2\left<\mu,\alpha_i\right>}{\left<\alpha_i,\alpha_i\right>}\alpha_i,\end{equation} and let $W=\left<s_1,...,s_n\right>$ denote the Weyl group. The standard presentation for the braid group $B$ is the generating set $\{T_w :w\in W\}$ subject to the relations $T_w T_{w^\prime} = T_{ww^\prime}$ for every $w,w^\prime\in W$ satisfying $\ell(w)+\ell(w^\prime)=\ell(ww^\prime)$, where $\ell$ is the length function on $W$. For each $i\in\{1,...,n\}$, we set $T_i:=T_{s_i}$. Thus, the braid group $B$ is generated by $T_1,...,T_n$. When $q$ is not a root of unity, $B$ acts via algebra automorphisms on ${\mathcal U}$ as follows:
\begin{align}
&T_i v_{\mu} =v_{s_i(\mu)},   \hspace{.6cm}  T_i  u_i^+ = -u_i^-v_{\alpha_i},\hspace{.6cm}   T_i  u_i^- = -v_{-\alpha_i}u_i^+,   &  \\\
&T_iu_j^+ = \sum_{r=0}^{-c_{ij}}\frac{(-q_i)^{-r}}{[-c_{ij}-r]_{q_i}![r]_{q_i}!}(u_i^+)^{-c_{ij}-r}u_j^+(u_i^+)^r ,   &   (i\neq j),  \\
&T_iu_j^- = \sum_{r=0}^{-c_{ij}}\frac{(-q_i)^r}{[-c_{ij}-r]_{q_i}![r]_{q_i}!}(u_i^-)^ru_j^-(u_i^-)^{-c_{ij}-r} ,   &    (i\neq j).
\end{align}
for all $i,j\in\{1,...,n\}$, $\mu\in Q$ \cite{L}.

Fix $w\in W$. For a reduced expression \begin{equation}w = s_{i_1}\cdots s_{i_t}\end{equation} define the roots \begin{equation}\beta_1 = \alpha_{i_1}, \beta_2=s_{i_1}\alpha_{i_2},...,\beta_t = s_{i_1}\cdots s_{i_{t-1}}\alpha_{i_t}\end{equation} and the root vectors \begin{equation}X_{\beta_1} = u_{i_1}^+, X_{\beta_2}=T_{s_{i_1}}u_{i_2}^+,...,X_{\beta_t}=T_{s_{i_1}}\cdots T_{s_{i_{t-1}}}u_{i_t}^+.\end{equation} Following \cite{DKP}, let ${\mathcal U}_q^w$ denote the subalgebra of ${\mathcal U}$ generated by the root vectors $X_{\beta_1},...,X_{\beta_t}$ (depends on the reduced expression).

When $k$ is algebraically closed of characteristic zero and $q$ is not  a root of unity, De Concini, Kac, and Procesi proved the following:

\begin{theorem} \cite[Proposition 2.2]{DKP}  If $(c_{ij})$ is a finite-type Cartan matrix, then the algebra ${\mathcal U}_q^w$ does not depend on the reduced expression for $w$. The algebra ${\mathcal U}_q^w$ has the PBW basis \[ X_{\beta_1}^{n_1}\cdots X_{\beta_t}^{n_t},\hspace{.4cm} n_1,...,n_t\in\mathbb{Z}_{\geq 0}.\] \end{theorem}

Beck later proved the analogous result for the case when $(c_{ij})$ is an affine Cartan matrix (with $q$ transcendental over $\mathbb{Q}$) \cite{B2}.

\section{A smash product of type \texorpdfstring{$D_{n+1}$}{D\_\{n+1\}}}\label{section smashD}

\subsection{The Algebras ${\mathcal U}_q(\mathfrak{so}_{2n+2})$, ${\mathcal U}_q^{w_n}$, and ${\mathcal O}_q(\mathfrak{o}k^{2n})$}

Fix an integer $n\geq 3$, and let $Q(D_{n+1})$ be the additive abelian subgroup of $\mathbb{R}^{n+1}$ consisting of the vectors having integer-valued coordinates $(a_1,...,a_{n+1})$ with the sum $\sum a_i$ being an even number. Let $\left<\,,\,\right>$ denote the restriction of the standard inner product on $\mathbb{R}^{n+1}$ (i.e. $\left<e_i,e_j\right>=\delta_{ij}$) to  $Q(D_{n+1})$. The group $Q(D_{n+1})$ is generated by the positive simple roots $\alpha_i = e_i-e_{i-1}$ for $2\leq i\leq n+1$ and $\alpha_1= e_1+e_2$.  For a positive simple root $\alpha_i$, let $s_i$ denote the corresponding simple reflection and let $W(D_{n+1})=\left<s_1,...,s_{n+1}\right>$ denote the associated Weyl group. The associated Cartan matrix $(c_{ij})$ is symmetric. Hence $d_1=\cdots =d_{n+1}=1$. Therefore, the parameters  $q_1,...,q_{n+1}$ are all equal to $q$. As usual, we put $\hat{q} = q-q^{-1}$.  Let ${\mathcal U}_q(\mathfrak{so}_{2n+2})$ denote the corresponding quantized universal enveloping algebra. We label the generators of  ${\mathcal U}_q(\mathfrak{so}_{2n+2})$  by $E_1,...,E_{n+1},F_1,...,F_{n+1}$ and $\{K_\mu:\mu\in Q(D_{n+1})\}$ and the defining relations are
\begin{align}
&\label{Uqrelns1}K_0=1,                                                                                                \hspace{.7cm}   K_\mu K_{\lambda}=K_{\mu+\lambda},                                   &  \\
&\label{Uqrelns2}K_{\mu} E_{i}=q^{\left<\mu,\alpha_i\right>}E_iK_{\mu},            \hspace{.7cm} K_{\mu} F_{i}=q^{-\left<\mu,\alpha_i\right>}F_iK_{\mu},       &  \\
&\label{Uqrelns3}E_iE_j=E_jE_i,                                                                                 \hspace{.7cm} F_iF_j=F_jF_i,                                                                              &\text{ ($\left<\alpha_i,\alpha_j\right>=0$ or $2$),}\\
&\label{Uqrelns4}E_i[E_i,E_j]=q[E_i,E_j]E_i,                                                              \hspace{.7cm} F_i[F_i,F_j]=q[F_i,F_j]F_i,                                                          &\text{ ($\left<\alpha_i,\alpha_j\right>=-1$),}   \\
&\label{Uqrelns5}E_iF_j\!=\!F_jE_i\!+\!\frac{\delta_{ij}}{\hat{q}}(K_{\alpha_i}\!\!-\!K_{\!-\alpha_i}),                                                                                                        &
\end{align}
for every $i,j\in\{1,...,n+1\}$ and $\mu,\lambda\in Q(D_{n+1})$. Here we use the $q^{-1}$-commutators, defined by \[ [u,v] := uv-q^{-1}vu \] for every $u,v\in{\mathcal U}_q(\mathfrak{so}_{2n+2})$.

Let $w_0$ denote the longest element of $W(D_{n+1})$ and let $w_0^L$ be the longest element of the parabolic subgroup $\left<s_1,...,s_{n}\right>\subseteq W(D_{n+1})$. Put $w_{n}=w_0^Lw_0$. We have a reduced expression \begin{equation}w_{n} = (s_{n+1}\cdots s_2s_1)(s_3\cdots s_{n}s_{n+1})\in W(D_{n+1})\end{equation} and root vectors \begin{equation}X_{e_{n+1}-e_{n}},X_{e_{n+1}-e_{n-1}},...,X_{e_{n+1}-e_1},X_{e_{n+1}+e_1}X_{e_{n+1}+e_2},...,X_{e_{n+1}+e_{n}}.\end{equation} For brevity, we put $x_i=X_{e_{n+1}-e_i}$ and $y_i=X_{e_{n+1}+e_i}$ for every $i\in\{1,...,n\}$. Let ${\mathcal U}_q^{w_n}$ denote the corresponding DeConcini-Kac-Procesi algebra.

The following can be found in \cite{GY}, Section 5.6.a.
\begin{theorem}The algebra ${\mathcal U}_q^{w_n}$ is isomorphic to the even-dimensional quantum Euclidean space ${\mathcal O}_q(\mathfrak{o}k^{2n})$.\end{theorem}
\begin{proof}We observe that the root vectors $x_1,...,x_{n},y_1,...,y_{n}$ of ${\mathcal U}_q^{w_n}$ can be written inductively as $x_{n}=E_{n+1}, y_1 = [x_2,E_1] $ and
\begin{align}\label{induction x}&x_i = [x_{i+1},E_{i+1}], \\
                    \label{induction y}&y_{i+1} = [y_i,E_{i+1}],
\end{align}
for all $1\leq i< n$. Using these identities, one can readily check that the root vectors satisfy the defining relations of ${\mathcal O}_q(\mathfrak{o}k^{2n})$ (c.f. \cite[Section 9.3.2]{KS97}),
\begin{align}
&\label{quantum euclidean relns1} x_ix_j=q^{-1}x_jx_i, \hspace{1cm}y_iy_j=qy_jy_i,    &(1\leq i<j\leq n),\\
&\label{quantum euclidean relns2} x_iy_j = q^{1-\delta_{ij}}y_jx_i +\delta_{ij}\hat{q}\sum_{r=1}^{i-1} (-q)^{i-r-1} x_ry_r,   &(i,j\in\{1,...,n\}).
\end{align}

Since ${\mathcal U}_q^{w_n}$ has a PBW basis of ordered monomials, Eqns. \ref{quantum euclidean relns1} and \ref{quantum euclidean relns2} are the defining relations. Hence, ${\mathcal U}_q^{w_n}\cong {\mathcal O}_q(\mathfrak{o}k^{2n})$.
\end{proof}

\subsection{${\mathcal U}_q^{w_n}$ as a left ${\mathcal U}_D^{\geq 0}$-module algebra}

Let ${\mathcal U}_D^{\geq 0}$ be the sub-Hopf algebra of ${\mathcal U}_q(\mathfrak{so}_{2n+2})$ generated by $E_1,...,E_{n+1}$, and $K_\mu$ for all $\mu\in Q(D_{n+1})$. We let $\pi:{\mathcal U}_D^{\geq 0}\to{\mathcal U}_D^{\geq 0}$ be the unique algebra map such that \begin{align}&\pi(E_{n+1})=0,\\  &\pi(E_i)=E_i&(i\leq n),\\   &\pi(K_\mu)=K_\mu&(\mu\in Q(D_{n+1})).\end{align} We define a function $\lambda : {\mathcal U}_D^{\geq 0}\otimes {\mathcal U}_q^{w_n}\to {\mathcal U}_D^{\geq 0}$ by  the following sequence of linear maps: \begin{equation}\label{definition of lambda}\xymatrix{\lambda : {\mathcal U}_D^{\geq 0}\otimes {\mathcal U}_q^{w_n}\ar[r]^{\hspace{0.5cm}incl.}&\Big({\mathcal U}_D^{\geq 0}\Big)^{\otimes 2}\ar[r]^{\pi\otimes id}& \Big({\mathcal U}_D^{\geq 0}\Big)^{\otimes 2}\ar[r]^{\hspace{.4cm}adjoint}&{\mathcal U}_D^{\geq 0}}\end{equation} and have the following:
\begin{theorem}\label{action is defined} For the function $\lambda$ above, we have $\text{Im}(\lambda)\subseteq {\mathcal U}_q^{w_n}$. In particular, $\lambda$ endows ${\mathcal U}_q^{w_n}$ with the structure of a left ${\mathcal U}_D^{\geq 0}$-module algebra.\end{theorem}
\begin{proof} For brevity, we set $u.v=\lambda(u\otimes v)$ for every $u\in{\mathcal U}_D^{\geq 0}$ and $v\in{\mathcal U}_q^{w_n}$. One can verify that
\begin{align}\label{action on x_r}E_j. x_r &=\begin{cases}-q(\delta_{1r}y_2+\delta_{2r}y_1), & (j=1), \\
                                                                     -q\delta_{jr}x_{r-1},                            & (j\neq 1),
                                                  \end{cases} \\
\label{action on y_r}E_j. y_r &=\begin{cases} 0,                                               & (j=n+1), \\
                                                                        -q\delta_{j,r+1}y_{r+1},                   &(j\neq n+1),
                                              \end{cases}
 \end{align}
for all $r\in\{1,...,n\}, j\in\{1,...,n+1\}$. Since ${\mathcal U}_D^{\geq 0}$ is a left ${\mathcal U}_D^{\geq 0}$-module algebra (with respect to the adjoint action), the equations \ref{action on x_r} and \ref{action on y_r} above, together with the fact that the $K_\mu$'s act diagonally on ${\mathcal U}_q^{w_n}$, prove the desired result.\end{proof} Using the action map $\lambda$, we form the smash product algebra ${\mathcal U}_q^{w_n}\#{\mathcal U}_D^{\geq 0}$ and define the following subalgebra \begin{equation}\label{definition of smashD} ({\mathcal U}_q^{w_n})^\# := \left<1\#u, u\#1\mid u\in{\mathcal U}_q^{w_n}\right>\subseteq{\mathcal U}_q^{w_n}\#{\mathcal U}_D^{\geq 0}.\end{equation}Loosely speaking, we can think of $({\mathcal U}_q^{w_n})^\#$ as being a smash product of ${\mathcal U}_q^{w_n}$ with itself.  Observe for example that $({\mathcal U}_q^{w_n})^\#$ is isomorphic as a vector space to ${\mathcal U}_q^{w_n}\otimes {\mathcal U}_q^{w_n}$.

\subsection{A Presentation of $\left({\mathcal U}_q^{w_n}\right)^\#$}

We will spend the rest of this section giving an explicit presentation for the algebra $\left({\mathcal U}_q^{w_n}\right)^\#$ because this will be necessary for proving the main result of Section \ref{section affD} (Thm. \ref{affineD=smash}).

The algebra $({\mathcal U}_q^{w_n})^\#$ is generated by $1\#x_i$, $1\#y_i$, $x_i\# 1$, $y_i\#1$ for $i\in\{1,...,n\}$. To compute the relations among these generators, we need comultiplication formulas for the root vectors $x_1,...,x_{n},y_1,...,y_{n}\in{\mathcal U}_q^{w_n}$. First, we must introduce the elements $\epsilon_{ij}, E_{r\downarrow s},E_{s\uparrow r}\in{\mathcal U}_D^{\geq 0}$ for every $i,j\in\{1,...,n\}$ and $r,s\in\{1,...,n+1\}$ with $r\geq s$. They are defined recursively via
\begin{align}
E_{r\downarrow s} &= \begin{cases} E_r, & (r=s),\\ [E_{r\downarrow s+1}, E_s], & (r\neq s),\end{cases}\\
E_{s\uparrow r} &= \begin{cases} E_s, & (r=s),\\ [E_{s\uparrow r-1}, E_r], & (r\neq s),\end{cases}
\\
\epsilon_{1j} &=\begin{cases} 0,& (j=1),\\ T_jT_{j-1}\cdots T_2E_1, & (j\neq 1),\end{cases}\\
\epsilon_{i+1,j} &=\begin{cases}[\epsilon_{ij},E_{i+1}],& (j\neq i,i+1),\\
q\epsilon_{i,i+1}E_{i+1}-q^{-1}E_{i+1}\epsilon_{i,i+1}, & (j=i+1),\\
\epsilon_{i,i+1}+q^{-1}(\epsilon_{ii}E_{i+1}-E_{i+1}\epsilon_{ii}), &(j=i).\end{cases}
\end{align}
We have the following:
\begin{lemma}\label{comult formulas} For every $i\in\{1,...,n\}$,
\begin{align}\label{comult x}\Delta (x_i) &= K_{-\deg(x_i)}\otimes x_i+x_i\otimes 1 +\hat{q}\sum_{j = i+1}^nE_{j\downarrow i+1}K_{-\deg(x_j)}\otimes x_j,\\
                       \label{comult y}\Delta (y_i) &= K_{-\deg(y_i)}\otimes y_i+y_i\otimes 1 \\
                                                                       &\nonumber\hspace{.5cm}+\hat{q}\left(\sum_{j=1}^n\epsilon_{ij}K_{-\deg(x_j)}\otimes x_j+\sum_{j=1}^{i-1}E_{j+1\uparrow i}K_{-\deg(y_j)}\otimes y_j\right).
\end{align}
\end{lemma}
\begin{proof}Use the induction formulas from  Eqns. \ref{induction x} and \ref{induction y} together with the comultiplication formula given in Equation \ref{comultiplication formula}.\end{proof}

From Eqns. \ref{action on x_r} and \ref{action on y_r} it follows that for all $i,j,r\in\{1,..., n\}$,
\begin{align}
\label{first}&E_{j\downarrow i+1}.x_r = -q \delta_{jr}x_i,  &  &E_{j\downarrow i+1}.y_r =  (-q)^{j-i}   \delta_{ir}y_j,\\
&E_{j+1\uparrow i}.x_r = (-q)^{i-j}\delta_{ir}x_j,  &   &E_{j+1\uparrow i}.y_r = -q \delta_{jr} y_i,\\
&\epsilon_{ij}.x_r = (-q)^{i+j-2}q^{\delta_{ij}}\delta_{ir}y_j -q \delta_{jr} y_i, &  \label{last}&\epsilon_{ij}.y_r = 0.
\end{align}
Using the identities \ref{first}-\ref{last} together with the comultiplication formulas, \ref{comult x}-\ref{comult y}, we compute the following ``cross-relations" in $({\mathcal U}_q^{w_n})^\#$.

\begin{proposition} For every $i,j\in\{1,...,n\}$,
\begin{align}
\label{smashrelns1}(1\# x_i)(x_j\# 1) &=\begin{cases} q^{-1}x_j\# x_i-q^{-1}\hat{q} x_i\# x_j, & i<j,\\
q^{-2} x_j\# x_i, & i=j,\\
q^{-1} x_j\# x_i, & i>j,
\end{cases}\\
\label{smashrelns2}(1\# y_i)(y_j\# 1) &=\begin{cases} q^{-1}y_j\# y_i-q^{-1}\hat{q} y_i\# y_j, & i>j,\\
q^{-2} y_j\# y_i, & i=j,\\
q^{-1} y_j\# y_i, & i<j,
\end{cases} \\
\label{smashrelns3}(1\# y_i)(x_j\# 1) &=  q^{-1+\delta_{ij}}x_j\#y_i-\hat{q}q^{-1}y_i\# x_j\\
&\nonumber\phantom{=}+\hat{q}q^{\!-\!1}\delta_{ij}\!\left(\!\sum_{m=1}^n\!(-q)^{i+m-2}y_m\# x_m
\!+\!\sum_{m=1}^{i-1}(-q)^{i-m}x_m\# y_m\!\right),\\
\label{smashrelns4}(1\# x_i)(y_j\# 1) &= q^{-1+\delta_{ij}}y_j\# x_i +\hat{q}q^{-1}\delta_{ij}\sum_{m=i+1}^n(-q)^{m-i}y_m\# x_m.
\end{align}
\end{proposition}

We have the following presentation for $({\mathcal U}_q^{w_n})^\#$:
\begin{theorem} The algebra $({\mathcal U}_q^{w_n})^\#$ is generated by $1\#x_i$, $x_i\#1$ for $1\leq i \leq n$, and its defining relations are Eqns. \ref{smashrelns1}-\ref{smashrelns4}  together with the relations
\begin{align}
\label{smashrelns5} (1\#x_i)(1\#x_j) &=q^{-1}(1\#x_j)(1\#x_i), &(1\leq i<j\leq n),\\
\label{smashrelns6} (1\#y_i)(1\#y_j) &=q(1\#y_j)(1\#y_i),    &(1\leq i<j\leq n),\\
\label{smashrelns7} (1\#x_i)(1\#y_j) &= q^{1-\delta_{ij}}(1\#y_j)(1\#x_i) & \\
\nonumber &\phantom{=}+\delta_{ij}\hat{q}\sum_{r=1}^{i-1} (-q)^{i-r-1} (1\#x_r)(1\#y_r),   &(i,j\in\{1,...,n\}),\\
\label{smashrelns8} (x_i\#1)(x_j\#1) &=q^{-1}(x_j\#1)(x_i\#1), &(1\leq i<j\leq n),\\
\label{smashrelns9} (y_i\#1)(y_j\#1) &=q(y_j\#1)(y_i\#1),    &(1\leq i<j\leq n),\\
\label{smashrelns10} (x_i\#1)(y_j\#1) &= q^{1-\delta_{ij}}(y_j\#1)(x_i\#1) & \\
\nonumber &\phantom{=}+\delta_{ij}\hat{q}\sum_{r=1}^{i-1} (-q)^{i-r-1} (x_r\#1)(y_r\#1),   &(i,j\in\{1,...,n\}).
\end{align}
\end{theorem}
\begin{proof} The generators $1\#x_1,...,1\#x_n$ generate a subalgebra isomorphic to ${\mathcal U}_q^{w_n}$, as do the generators $x_1\#1,...,x_n\#1$, giving us the relations \ref{smashrelns5}-\ref{smashrelns10}. The universal property of smash products (for example, see \cite[Section 1.8]{HS}) and the PBW basis of De Concini-Kac-Procesi algebras imply that the cross relations of \ref{smashrelns1}-\ref{smashrelns4} together with the above relations are a presentation of  $({\mathcal U}_q^{w_n})^\#$.
\end{proof}

\section{The quantum affine algebra \texorpdfstring{${\mathcal U}_q^{\widehat{w}_n}$}{U\_q\{w\_n\}}}\label{section affD}

Let $Q(\widehat{D}_{n+1})=Q(D_{n+1})\oplus\mathbb{Z}$ denote the root lattice of type $\widehat{D}_{n+1}$. As an abelian group, $Q(\widehat{D}_{n+1})$ is generated additively by the positive simple roots $\alpha_0:=-e_{n+1}-e_n+1$, $\alpha_1:=e_1+e_2$, and $\alpha_i:=e_i-e_{i-1}$ for $2\leq i\leq n+1$. We extend the bilinear form $\left<\,\,\, ,\,\,\,\right>$ on $Q(D_{n+1})$ to  $Q(\widehat{D}_{n+1})$  by setting $1\in Q(\widehat{D}_{n+1})$ to be isotropic. As before, let $s_i$ denote the corresponding simple reflection $s_i : Q(\widehat{D}_{n+1})\to Q(\widehat{D}_{n+1})$, for $0\leq i\leq n+1$, and $W(\widehat{D}_{n+1})=\left<s_0,...,s_{n+1}\right>$ is the Weyl group. The corresponding quantized enveloping algebra ${\mathcal U}_q(\widehat{\mathfrak{so}}_{2n+2})$ is generated by $E_0,...,E_{n+1},F_0,...,F_{n+1}$ and $\{K_\mu:\mu\in Q(\widehat{D}_{n+1})\}$ and has defining relations
\begin{align}
&\label{Uqrelns1_aff}K_0=1,                                                                                                \hspace{.7cm}   K_\mu K_{\lambda}=K_{\mu+\lambda},                                   &  \\
&\label{Uqrelns2_aff}K_{\mu} E_{i}=q^{\left<\mu,\alpha_i\right>}E_iK_{\mu},            \hspace{.7cm} K_{\mu} F_{i}=q^{-\left<\mu,\alpha_i\right>}F_iK_{\mu},       &  \\
&\label{Uqrelns3_aff}E_iE_j=E_jE_i,                                                                                 \hspace{.7cm} F_iF_j=F_jF_i,                                                                              &\text{ ($\left<\alpha_i,\alpha_j\right>=0$ or $2$),}\\
&\label{Uqrelns4_aff}E_i[E_i,E_j]=q[E_i,E_j]E_i,                                                              \hspace{.7cm} F_i[F_i,F_j]=q[F_i,F_j]F_i,                                                          &\text{ ($\left<\alpha_i,\alpha_j\right>=-1$),}   \\
&\label{Uqrelns5_aff}E_iF_j\!=\!F_jE_i\!+\!\frac{\delta_{ij}}{\hat{q}}(K_{\alpha_i}\!\!-\!K_{\!-\alpha_i}),                                                                                                        &
\end{align}
for every $i,j\in\{0,...,n+1\}$ and $\mu,\lambda\in Q(\widehat{D}_{n+1})$ (c.f. Eqns. \ref{Uqrelns1}-\ref{Uqrelns5}).

Let $\widehat{w}_n\in W(\widehat{D}_{n+1})$ be the Weyl group element given by \begin{equation}\widehat{w}_n: v + r \mapsto v + r+2a_{n+1} \end{equation} for every $v=\sum_{i=1}^{n+1} a_ie_i\in Q(D_{n+1})$ and $r\in\mathbb{Z}$. We have the reduced expression \begin{equation}\widehat{w}_n:=(s_{n+1}\cdots s_1)(s_3\cdots s_{n+1})s_0(s_n\cdots s_3)(s_1\cdots s_n)s_0\in W(\widehat{D}_{n+1}).\end{equation}

We let $\widehat{B}_{\mathfrak{so}_{2n+2}}=\left<T_0,...,T_{n+1}\right>$  denote the corresponding braid group of $\widehat{\mathfrak{so}}_{2n+2}$ and label the corresponding ordered root vectors for ${\mathcal U}_q^{\widehat{w}_n}$ by \begin{equation}X_n,...,X_1,Y_1,...,Y_n,\overline{X}_n,...,\overline{X}_1,\overline{Y}_1,...,\overline{Y}_n.\end{equation}

One can readily verify the following lemmas.

\begin{lemma}\label{lemma recursion formulas affD} We have the following recursion formulas in the algebra ${\mathcal U}_q^{\widehat w_n}$:
\begin{align}
&X_n                           =E_{n+1},                                      &X_i                             &=[X_{i+1},E_{i+1}],                      &   \text{($i\neq n$),}&        \\
&Y_1                           =[X_2,E_1],                                   &Y_i                            &=[Y_{i-1},E_i],                               &    \text{($i\neq 1$),}&            \\
&\overline{X}_n         =[Y_{n-1},T_{n+1}T_nE_0],       &\overline{X}_i           &=[\overline{X}_{i+1},E_{i+1}],      &   \text{($i\neq n$),}&        \\
&\overline{Y}_1         =[\overline{X}_2,E_1],                &\overline{Y}_i           &=[\overline{Y}_{i-1},E_i],              &    \text{($i\neq 1$),}&           \\
& Y_2                          =[X_1,E_1],                                  &\overline{Y}_2          &=[\overline{X}_1,E_1].                 &     &
\end{align}
\end{lemma}

\begin{lemma}\label{lemma braid group action on affD root vectors} For all $i,j\in\{1,...,n\}$, we have the following:
\begin{align}
T_i.X_j&=\begin{cases}  [E_i,X_j],  &    (i=j \text{ or } (i,j)=(1,2)), \\
                                          X_{j+1},                  & (i=j+1),\\
                                          X_j,                         & \text{otherwise},
                                          \end{cases}
\\
T_i.Y_j&=\begin{cases}  Y_{j-1},  &    (i=j\text{ and } i\neq 1), \\
                                          X_{3-j},   &  (i=1,j\in\{1,2\}),\\
                                          [E_{j+1},Y_j],                  & (i=j+1),\\
                                          Y_j,                         & \text{otherwise},
                                          \end{cases}
\\
T_i.\overline{X}_j&=\begin{cases}  [E_i,\overline{X}_j],  &    (i=j \text{ or } (i,j)=(1,2)), \\
                                          \overline{X}_{j+1},                  & (i=j+1),\\
                                          \overline{X}_j,                         & \text{otherwise},
                                          \end{cases}
\\
T_i.\overline{Y}_j&=\begin{cases}  \overline{Y}_{j-1},  &    (i=j\text{ and } i\neq 1), \\
                                          \overline{X}_{3-j},   &  (i=1,j\in\{1,2\}),\\
                                          [E_{j+1},\overline{Y}_j],                  & (i=j+1),\\
                                          \overline{Y}_j,                         &\text{otherwise}.
                                          \end{cases}\end{align}
\end{lemma}

With the help of Lemmas \ref{lemma recursion formulas affD} and \ref{lemma braid group action on affD root vectors}, we prove the following.

\begin{proposition} The defining relations for the algebra ${\mathcal U}_q^{\widehat w_n}$ are
\begin{align}
&\label{affinerelns1}X_iX_j=q^{-1}X_jX_i,\hspace{1 cm}Y_jY_i=q^{-1}Y_iY_j,&(i<j),\\
&\label{affinerelns2}\overline X_i\overline X_j=q^{-1}\overline X_j\overline X_i,\hspace{1 cm}\overline Y_j\overline Y_i=q^{-1}\overline Y_i\overline Y_j,&(i<j),\\
&\label{affinerelns3}Y_jX_i = q^{\delta_{ij}-1}X_iY_j-\delta_{ij}\hat{q}\sum_{r=1}^{i-1}(-q)^{i-r-1}X_{r}Y_{r},&\\
&\label{affinerelns4}\overline{Y}_j\overline{X}_i = q^{\delta_{ij}-1}\overline{X}_i\overline{Y}_j-\delta_{ij}\hat{q}\sum_{r=1}^{i-1}(-q)^{i-r-1}\overline{X}_{r}\overline{Y}_{r},&\\
&\label{affinerelns5}\overline{X}_iX_i=q^{-2}X_i\overline{X}_i,\hspace{1 cm}\overline{Y}_iY_i =q^{-2}Y_i\overline{Y}_i,  &\\
&\label{affinerelns6}\overline{X}_jX_i=q^{-1}X_i\overline{X}_j, \hspace{1 cm}\overline Y_iY_j=q^{-1}Y_j\overline Y_i,&(i<j),\\
&\label{affinerelns7}\overline{X}_iX_j=q^{-1}X_j\overline{X}_i-q^{-1}\hat{q}X_i\overline{X}_j,\hspace{.3 cm} \overline{Y}_jY_i=q^{-1}Y_i\overline{Y}_j-q^{-1}\hat{q}Y_j\overline{Y}_i,&(i<j),\\
&\label{affinerelns8}\overline X_iY_j=q^{-1+\delta_{ij}}Y_j\overline X_i+\widehat qq^{-1}\delta_{ij}\sum_{m=i+1}^n(-q)^{m-i}Y_m\overline X_m,  &\\
&\label{affinerelns9}\overline Y_iX_j      =q^{-1+\delta_{ij}}X_j\overline Y_i-\widehat qq^{-1} Y_i\overline X_j    & \\
&\nonumber \hspace{1.2cm}+\widehat qq^{-1}\delta_{ij}\!\!\left[\sum_{m=1}^n(-q)^{i+m\!-\!2}Y_m\overline X_m \!+\!\sum_{m=1}^{i-1}(-q)^{i-m}X_m\overline Y_m\!\right],   &
\end{align}
for $i,j\in\{1,...,n\}$.
\end{proposition}
\begin{proof} The first $2n$ letters in the reduced expression for $\widehat{w}_n$ coincide with $w_n$, as do the last $2n$ letters. This gives us the relations \ref{affinerelns1}-\ref{affinerelns4}. Using Lemmas \ref{lemma recursion formulas affD} and \ref{lemma braid group action on affD root vectors}, one can prove inductively that the remaining relations hold. To illustrate how to obtain the identities in Eqn. \ref{affinerelns5} for example, one can first verify the base cases, $\overline{X}_1X_1 = q^{-2}X_1\overline{X}_1$ and $\overline{Y}_nY_n=q^{-2}Y_n\overline{Y}_n$, and then apply appropriate braid group automorphisms (refer to Lemma \ref{lemma braid group action on affD root vectors}) to both sides of the equations. Since ${\mathcal U}_q^{\widehat w_n}$ has a PBW basis of ordered monomials, Eqns. \ref{affinerelns1}-\ref{affinerelns9} are the defining relations. \end{proof}

By comparing Eqns. \ref{smashrelns1}-\ref{smashrelns10} with Eqns. \ref{affinerelns1}-\ref{affinerelns9}, we observe the following theorem.
\begin{theorem}\label{affineD=smash}As $k$-algebras, ${\mathcal U}_q^{\widehat w_n}\cong ({\mathcal U}_q^{w_n})^\#$ via the isomorphism $$\begin{matrix} &X_i\mapsto (x_i\# 1),&Y_i\mapsto (y_i\#1),\\ &\overline X_i\mapsto (1\# x_i),&\overline Y_i\mapsto (1\#y_i),\\ \end{matrix}\ \ \  \text{ for }i=1,\dots,n.$$ \end{theorem}

\section{The FRT-Construction and the algebra \texorpdfstring{$\mathbb{X}_{n,q}$}{X\_\{n,q\}}}\label{section X algebra}
We will briefly review the Faddeev-Reshetikhin-Takhtajan (FRT) construction of \cite{FRT88} (see \cite[Section 7.2]{CP} for more details). We let $V$ be a $k$-module with basis $\{v_1,\dots,v_N\}$.  For a linear map $R\in\text{End}_k(V\otimes V)$, we write \begin{equation}R(v_i\otimes v_j)=\sum_{s,t}R_{ij}^{st}v_s\otimes v_t\text{ for all }1\leq i,j<N,\end{equation} with all $R_{ij}^{st}\in k$. The \emph{FRT algebra} $\A(R)$ associated to $R$ is the $k$-algebra presented by generators $X_{ij}$ for $1\leq i,j\leq N$ and has the defining relations \begin{equation}\label{FRT relations}\sum_{s,t}R_{st}^{ji}X_{sl}X_{tm}=\sum_{s,t}R_{lm}^{ts}X_{is}X_{jt}\end{equation} for every $i,j,l,m\in\{1,...,N\}$. Up to algebra isomorphism, $\A(R)$ is independent of the chosen basis of $V$.

Let us specialize now to the case when $N=2n$. Following \cite[Section 8.4.2]{KS97}, for each $i,j\in\{1,\dots, 2n\}$, let $E_{ij}$ denote the linear map on $V$ defined by  $E_{ij}.v_\ell = \delta_{j\ell} v_i$. Let $i^\prime :=2n+1-i$, and let
\begin{align}
\label{RmatrixD}R_{D_n}&=q\sum_{i:i\neq i^\prime}(E_{ii}\otimes E_{ii})+\sum_{i,j:i\neq j,j'}(E_{ii}\otimes E_{jj})+q^{-1}\sum_{i:i\neq i'}(E_{i'i'}\otimes E_{ii})\\
\nonumber &\phantom{=}+\hat q\left(\sum_{i,j:i>j}(E_{ij}\otimes E_{ji})-\sum_{i,j:i>j}q^{\rho_i-\rho_j}(E_{ij}\otimes E_{i'j'})\right),
\end{align}
where $(\rho_1,\rho_2,\dots,\rho_{2n})$ is the $2n$-tuple $(n-1,n-2,\dots,1,0,0,-1,\dots,-n+1).$

We define an algebra $\X_{n,q}$ presented by generators $X_{ij}$ with $i\in\{1,2\}$, $j\in\{1,...,2n\}$, and having the defining relations
\begin{align}
\label{r1}X_{rt}X_{rs}&=q^{-1} X_{rs}X_{rt}&(r\in\{1,2\},s<t,\ t\neq s'),\\
X_{rs'}X_{rs}&=X_{rs}X_{rs'}+\hat q\sum_{l=s+1}^nq^{l-s-1}X_{rl}X_{rl'}&(r\in\{1,2\},s<s'),\\
X_{2s}X_{1s}&=q^{-1} X_{1s}X_{2s},\\
X_{2s}X_{1t}&=X_{1t}X_{2s}&(s<t,\ t\neq s'),\\
X_{2t}X_{1s}&=X_{1s}X_{2t}-\hat qX_{1t}X_{2s}&(s<t,\ t\neq s'),\\
X_{2s}X_{1s'}&=qX_{1s'}X_{2s}+\hat q\sum_{l=1}^{s-1}q^{s-l}X_{1l'}X_{2l}&(s<s'),\\
\label{r2}X_{2s'}X_{1s}&=qX_{1s}X_{2s'}+\hat q\sum_{l=s+1}^nq^{l-s}X_{1l}X_{2l'},\\
&\nonumber\phantom{=}+\hat{q}q^{-1}\sum_{l=1}^nq^{l'-s}X_{1l'}X_{2l}-\hat qX_{1s'}X_{2s}&(s<s').
\end{align}

We label the canonical generators of $\A(R_{D_n})$ by $Y_{ij}$ for $i,j=1,\dots,2n$, and let $T_{2,n}$ be the subalgebra of $\A(R_{D_n})$ generated by $\{Y_{ij}:1\leq i\leq 2,\ 1\leq j\leq 2n\}$.

\begin{proposition}\label{surjective homomorphism} There is a surjective algebra homomorphism $\X_{n,q}\to T_{2,n}$ with kernel $\left<\Omega_1,\Omega_2,\Upsilon\right>$, where
\begin{align}
\label{extrarelns}\Omega_1:=\sum_{r=1}^nq^{\rho_{r'}}X_{1,r}X_{1,r'}, & &
\Omega_2:=\sum_{r=1}^nq^{\rho_{r'}}X_{2,r}X_{2,r'}, & &
\Upsilon:=\sum_{r=1}^{2n}q^{\rho_r}X_{1,r'}X_{2,r}.
\end{align}
\end{proposition}

\begin{proof}Using the FRT construction (see Equation \ref{FRT relations} and \ref{RmatrixD}), one can readily compute the defining relations for the  algebra $\A(R_{D_n})$ and see that they line up appropriately with Equation \ref{r1}-\ref{r2} together with $\Omega_1=\Omega_2=\Upsilon=0$.\end{proof}

Notice that the definition of ${\mathbb X}_{n,q}$ makes sense when $n=2$, and Proposition \ref{surjective homomorphism} holds in this case as well. However, the rest of the results of this paper require $n\geq 3$.

Following \cite{AST91}, we recall the details on twisting algebras by cocycles. Let $M$ be an additive abelian group and $c:M\times M\to k^\times$ a $2$-cocycle of $M$. If $\Lambda$ is a $k$-algebra graded by $M$, we can twist $\Lambda$ by $c$ to obtain a new $M$-graded $k$-algebra $\Lambda'$ that is canonically isomorphic to $\Lambda$ as a $k$-module via $x\leftrightarrow x'$.  Multiplication of homogeneous elements in $\Lambda'$ is given by $$x'y'=c(\deg(x),\deg(y))(xy)'.$$

For our purposes, we will let $\beta:Q(\widehat{D}_{n+1})\times Q(\widehat{D}_{n+1})\to k^\times$ be the bicharacter (hence, also a $2$-cocycle) defined by \begin{equation}\label{definition of beta}\beta(\alpha_i,\alpha_j)=\begin{cases} q & (i,j)=(0,n+1),\\ 1 &(i,j)\neq(0,n+1),\end{cases}\end{equation} and have the following:
\begin{theorem}\label{twisting} The $\beta$-twisted algebra $\left({\mathcal U}_q^{\widehat{w}_n}\right)^\prime$ is isomorphic to $\X_{n,q}$.\end{theorem}
\begin{proof} We label the corresponding generators of $\left({\mathcal U}_q^{\widehat w_n}\right)^\prime$ by \begin{equation}X_n^\prime,...,X_1^\prime,Y_1^\prime,...,X_n^\prime,\overline{X}_n^\prime,...,\overline{X}_1^\prime,\overline{Y}_1^\prime,...,\overline{Y}_n^\prime.\end{equation} By comparing Eqns. \ref{affinerelns1}-\ref{affinerelns9} and \ref{r1}-\ref{r2}, we observe that the algebra map $\left({\mathcal U}_q^{\widehat w_n}\right)^\prime\to\mathbb{X}_{n,q}$ defined by
\begin{align}
 X_i^\prime     &\mapsto (-1)^{n+1-i}X_{1,n+1-i},     &   &\overline{X}_i^\prime\mapsto (-1)^{n+1-i}X_{2,n+1-i},\\
 Y_i^\prime    &\mapsto X_{1,n+i}, &     &\overline{Y}_i^\prime\mapsto X_{2,n+i},
\end{align}
for every $i\in\{1,...,n\}$, is an isomorphism.
\end{proof}

From this, we deduce the following:
\begin{theorem} The algebra $\X_{n,q}$ is an iterated Ore extension over $k$,\[ \X_{n,q} = k[X_{11}][X_{12};\tau_{12},\delta_{12}]\cdots [X_{1n};\tau_{1n},\delta_{1n}][X_{21};\tau_{21},\delta_{21}]\cdots [X_{2n};\tau_{2n},\delta_{2n}].\]\end{theorem}
\begin{proof} It suffices to check that ordered monomials are linearly independent. From Theorem \ref{twisting}, we have a canonical vector space isomorphism ${\mathcal U}_q^{\widehat w_n}\to\X_{n,q}$ that preserves the ordered generating sets. Since ${\mathcal U}_q^{\widehat w_n}$ has a basis of ordered monomials, $\X_{n,q}$ does as well.\end{proof}

\section{A type \texorpdfstring{$A_m$}{A\_m} analogue}\label{section A analogue}

\subsection{The algebras ${\mathcal U}_q(\mathfrak{sl}_{m+1})$, ${\mathcal U}_q^{c_m}$, and ${\mathcal O}_q(k^m)$.}

Fix an integer $m> 1$. Let $Q(A_m)$ denote the abelian subgroup of $\mathbb{R}_{m+1}$  consisting of integral $(m+1)$-tuples $(a_1,...,a_{m+1})$ with the sum $\sum a_i$ equalling $0$. As a group, $Q(A_m)$ is generated by $\alpha_i:=e_i-e_{i+1}$ for $i\in\{1,...,m\}$. Let $W(A_m)$ and $B_{\mathfrak{sl}_{m+1}}$ denote the corresponding Weyl group and braid group, respectively. Let ${\mathcal U}_q(\mathfrak{sl}_{m+1})$ denote the corresponding quantum enveloping algebra, and let  ${\mathcal U}_A^{\geq 0}$ be the positive  Borel subalgebra of ${\mathcal U}_q(\mathfrak{sl}_{m+1})$. We consider the Coxeter element \begin{equation} c_m = s_1\cdots s_m\in W(A_m) \end{equation} and the associated De Concini-Kac-Procesi algebra ${\mathcal U}_q^{c_m}$. We label the root vectors in ${\mathcal U}_q^{c_m}$ by \begin{equation} z_1:= X_{e_1-e_2},\ z_2:=X_{e_1-e_3},...,z_m:=X_{e_1-e_{m+1}} \end{equation} and have the following
\begin{proposition} The root vectors $z_1,...,z_m$ satisfy the relations \begin{equation}\label{c_m relns} z_iz_j = qz_jz_i \end{equation} for all $i,j\in\{1,...,m\}$ with $i<j$. \end{proposition}

Since ${\mathcal U}_q^{c_m}$ has a PBW basis of ordered monomials, the relations of Eqn. \ref{c_m relns} are the defining relations for ${\mathcal U}_q^{c_m}$. In particular, we have the following well-known result (c.f. for example \cite{MC}):
\begin{corollary} The algebra ${\mathcal U}_q^{c_m}$ is isomorphic to the algebra of quantum affine space ${\mathcal O}_q(k^m)$. \end{corollary}

Denote by $\pi_A:{\mathcal U}_A^{\geq 0}\to{\mathcal U}_A^{\geq 0}$ the unique algebra map such that \begin{align}&\pi(E_{1})=0,\\ &\pi(E_i)=E_i&(1<i\leq m),\\ &\pi(K_\mu)=K_\mu&(\mu\in Q(A_m)).\end{align}  Let $\lambda_A : {\mathcal U}_A^{\geq 0}\otimes {\mathcal U}_q^{c_m}\to {\mathcal U}_A^{\geq 0}$ be defined by the following sequence of linear maps:\begin{equation}\label{definition of lambda_A}\xymatrix{\lambda_A : {\mathcal U}_A^{\geq 0}\otimes {\mathcal U}_q^{c_m}\ar[r]^{\hspace{0.5cm}incl.}&\Big({\mathcal U}_A^{\geq 0}\Big)^{\otimes 2}\ar[r]^{\pi_A\otimes id}& \Big({\mathcal U}_A^{\geq 0}\Big)^{\otimes 2}\ar[r]^{\hspace{.4cm}adjoint}&{\mathcal U}_A^{\geq 0}}.\end{equation}

The identities in Equation~\ref{action on x_r}  imply the following
\begin{corollary}The linear map $\lambda_A$ satisfies $\text{Im}(\lambda_A)\subseteq{\mathcal U}_q^{c_{m}}$. In particular, $\lambda_A$ endows the algebra ${\mathcal U}_q^{c_m}$ with the structure of a left ${\mathcal U}_A^{\geq 0}$-module algebra.\end{corollary} As before (see \ref{definition of smashD}), we use the action map $\lambda_A$ to construct the smash product ${\mathcal U}_q^{c_m}\#{\mathcal U}_A^{\geq 0}$ and let $({\mathcal U}_q^{c_m})^\#$ denote the subalgebra \begin{equation}({\mathcal U}_q^{c_m})^\#:= \left< 1\# u,u\#1\mid u\in{\mathcal U}_q^{c_m}\right>\subseteq {\mathcal U}_q^{c_m}\#{\mathcal U}_A^{\geq 0}.\end{equation}

\subsection{The quantum affine algebra ${\mathcal U}_q^{\widehat{c}_m}$}\label{section affA}

Let $Q(\widehat{A}_m)=Q(A_m)\oplus\mathbb{Z}$ denote the root lattice of type $\widehat{A}_m$. As an abelian group, $Q(\widehat{A}_m)$ is generated additively by the positive simple roots $\alpha_0:=e_m-e_1+1$, and $\alpha_i:=e_i-e_{i+1}$ for $i\in\{1,...,m\}$. We extend the inner product $\left<\,\,\, ,\,\,\,\right>$ on $Q(A_m)$ to an inner product on $Q(\widehat{A}_m)$  by setting $1\in Q(\widehat{A}_m)$ to be isotropic. We let $s_i$ denote the corresponding simple reflection $s_i : Q(\widehat{A}_m)\to Q(\widehat{A}_m)$, for $0\leq i\leq m$, and let $W(\widehat{A}_m)=\left<s_0,...,s_m\right>$ denote the corresponding affine Weyl group. We let ${\mathcal U}_q(\widehat{\mathfrak{sl}}_{m+1})$ denote the corresponding quantized enveloping algebra.

We set \begin{equation} \widehat{c}_m : = (s_1\cdots s_m)(s_0s_1\cdots s_{m-1})\in W(\widehat{A}_m) \end{equation} and note the following analogue of Theorem \ref{affineD=smash}.
\begin{theorem}\label{affineA=smash} As $k$-algebras, ${\mathcal U}_q^{\widehat{c}_m}\cong ({\mathcal U}_q^{c_m})^\#$.\end{theorem}
\begin{proof}Compute. One can use an analogous isomorphism of Thm. \ref{affineD=smash}.\end{proof}

Now let $V$ be a $k$-module with basis $\{v_1,...,v_m\}$, and for all $i,j,\ell\in\{1,...,m\}$, define linear maps $e_{ij}$ by the rule $e_{ij}.v_{\ell}=\delta_{j\ell}v_i$.

Set \begin{equation}\label{RmatrixA}R_{A_{m-1}}=q\sum_{i=1}^m(e_{ii}\otimes e_{ii})+\sum_{i\neq j}(e_{ii}\otimes e_{jj})+\hat q\sum_{i>j}(e_{ij}\otimes e_{ji}).\\ \end{equation} This is the standard $R$-matrix of type $A_{m-1}$ (see \cite[Section 8.4.2]{KS97}).

The algebra of $m\times m$ quantum matrices, denoted $\O_q(M_m(k))$, is the algebra $\A(R_{A_{m-1}})$ and was defined in \cite{FRT88}.  More generally, one considers $\ell\times p$ quantum matrices, denoted $\O_q(M_{\ell,p}(k))$, by looking at appropriate subalgebras of square quantum matrices.

We let $\gamma : Q(\widehat{A}_m)\times Q(\widehat{A}_m)\to k^\times$ be the bicharacter defined by \begin{equation}\gamma (\alpha_i,\alpha_j)   =   \begin{cases}      q,    &   (i,j)=(0,1),\\ 1,  &(i,j)\neq (0,1).   \end{cases}\end{equation} and have the following analogue of Thm. \ref{twisting}.
\begin{theorem}\label{twisting_A} Twisting ${\mathcal U}_q^{\widehat{c}_m}$ by $\gamma$ yields an algebra isomorphic to ${\mathcal O}_q(M_{2,m})$. \end{theorem}
\begin{proof}Compute (c.f. Thm. \ref{twisting}).\end{proof}

Theorem \ref{twisting_A}, together with Proposition \ref{surjective homomorphism}, allows us to view $\X_{n,q}$ as an orthogonal analogue of $2\times n$ quantum matrices.  The key distinction is that $\O_q(M_{2,n}(k))$ is a subalgebra of $\A(R_{A_{n-1}})$, whereas $\X_{n,q}$ is a parent of the analogous subalgebra $T_{2,n}\subseteq\A(R_{D_n})$.

\end{document}